\documentstyle[txmac,a4,
twoside,%
nocaphead,%
epsf,%
varthm,%
rotate,%
myrot,%
mypic,times,mathptm]{article}

\advance\oddsidemargin by -0.7cm
\advance\evensidemargin by -0.7cm
\advance\textwidth by 1.4cm

\def\mynewtheo#1#2{%
\newtheorem{@#1}{#2}[section]%
\newenvironment{#1}{\begin{@#1}\rm}{\end{@#1}}}

\mynewtheo{lem}{Lemma}
\mynewtheo{exer}{Exercise}
\mynewtheo{theo}{Theorem}
\mynewtheo{rem}{Remark}
\mynewtheo{defi}{Definition}
\mynewtheo{conj}{Conjecture}
\mynewtheo{corr}{Corollary}
\mynewtheo{prop}{Proposition}
\mynewtheo{question}{Question}
\mynewtheo{exam}{Example}

\newenvironment{eqn}{\begin{equation}}{\end{equation}}

\newenvironment{theorem}{\begin{theo}}{\end{theo}}
\newenvironment{conjecture}{\begin{conj}}{\end{conj}}
\newenvironment{corollary}{\begin{corr}}{\end{corr}}
\newenvironment{proposition}{\begin{prop}}{\end{prop}}

\begin{document}

\makeatletter

\newenvironment{myeqn*}[1]{\begingroup\def\@eqnnum{\reset@font\rm#1}%
\xdef\@tempk{\arabic{equation}}\begin{equation}\edef\@currentlabel{#1}}
{\end{equation}\endgroup\setcounter{equation}{\@tempk}\ignorespaces}

\newenvironment{myeqn}[1]{\begingroup\let\eq@num\@eqnnum
\def\@eqnnum{\bgroup\let\r@fn\normalcolor 
\def\normalcolor####1(####2){\r@fn####1#1}%
\eq@num\egroup}%
\xdef\@tempk{\arabic{equation}}\begin{equation}\edef\@currentlabel{#1}}
{\end{equation}\endgroup\setcounter{equation}{\@tempk}\ignorespaces}

\newcount\case@cnt
\newenvironment{caselist}{\case@cnt0\relax}{}

\newcommand{\mybin}[2]{\text{$\Bigl(\begin{array}{@{}c@{}}#1\\#2%
\end{array}\Bigr)$}}
\def\overtwo#1{\mbox{\small$\mybin{#1}{2}$}}
\newcommand{\mybr}[2]{\text{$\Bigl\lfloor\mbox{%
\small$\displaystyle\frac{#1}{#2}$}\Bigr\rfloor$}}
\def\mybrtwo#1{\mbox{\mybr{#1}{2}}}

\def\case{\advance\case@cnt by 1\relax{\bf\relax
\the\case@cnt\relax. Case.\ \ignorespaces}%
\edef\@currentlabel{\the\case@cnt}}

\def\myfrac#1#2{\raisebox{0.2em}{\small$#1$}\!/\!\raisebox{-0.2em}{\small$#2$}}

\def\myeqnlabel{\bgroup\@ifnextchar[{\@maketheeq}{\immediate
\stepcounter{equation}\@myeqnlabel}}

\def\@maketheeq[#1]{\def\theequation{#1}\@myeqnlabel}

\def\@myeqnlabel#1{%
{\edef\@currentlabel{\theequation}
\label{#1}\enspace\eqref{#1}}\egroup}

\def\epsfs#1#2{{\epsfxsize#1\relax\epsffile{#2.eps}}}




\author{A. Stoimenow\footnotemark[1]\\[2mm]
\small Ludwig-Maximilians University Munich, Mathematics\\
\small Institute, Theresienstra\ss e 39, 80333 M\"unchen, Germany,\\
\small e-mail: {\tt stoimeno@informatik.hu-berlin.de},\\
\small WWW\footnotemark[2]\enspace: {\hbox{\tt http://www.informatik.hu-berlin.de/%
\raisebox{-0.8ex}{\tt\~{}}stoimeno}}
}

{\def\thefootnote{\fnsymbol{footnote}}
\footnotetext[1]{Supported by a DFG postdoc grant.}
\footnotetext[2]{I am sometimes criticized that my papers are not
available. All my papers, including those referenced here, are
available on my webpage or by sending me an inquiry to the email
address I specified above. Hence, \em{please}, do not complain about
some paper being non-available before trying these two options.}
}

\title{\large\bf \uppercase{%
Vassiliev invariants and rational knots}\\[2mm]
\uppercase{of unknotting number one}}

\date{\large Current version: \today\ \ \ First version:
\makedate{28}{7}{1997}}


\maketitle

\makeatletter

\def\ex{\exists\,}
\def\fa{\forall\,}
\let\nb\nabla
\let\reference\ref
\let\ay\asymp
\let\pa\partial
\let\ap\alpha
\let\bt\beta
\let\zt\zeta
\let\Gm\Gamma
\let\gm\gamma
\let\de\delta
\let\dl\delta
\let\Dl\Delta
\let\eps\epsilon
\let\lm\lambda
\let\Lm\Lambda
\let\sg\sigma
\let\vp\varphi
\let\om\omega
\let\es\enspace

\let\sm\setminus
\let\tl\tilde
\def\dt{\det}
\def\sgn{\mathop {\operator@font sgn}}
\def\rk{\mathop {\operator@font rank}}
\def\ncap{\not\mathrel{\cap}}
\def\cf{\text{\rm cf}\,}
\def\Ra{\Rightarrow}
\def\lra{\longrightarrow}
\def\llra{\longleftrightarrow}
\def\ol{\overline}
\def\so{\Rightarrow}
\def\So{\Longrightarrow}
\let\ds\displaystyle
\def\lz{\linebreak[0]\verb}

\long\def\@makecaption#1#2{%
   \vskip 10pt
   {\let\label\@gobble
   \let\ignorespaces\@empty
   \xdef\@tempt{#2}%
   }%
   \ea\@ifempty\ea{\@tempt}{%
   \setbox\@tempboxa\hbox{%
      \fignr#1#2}%
      }{%
   \setbox\@tempboxa\hbox{%
      {\fignr#1:}\capt\ #2}%
      }%
   \ifdim \wd\@tempboxa >\captionwidth {%
      \rightskip=\@captionmargin\leftskip=\@captionmargin
      \unhbox\@tempboxa\par}%
   \else
      \hbox to\captionwidth{\hfil\box\@tempboxa\hfil}%
   \fi}%
\def\fignr{\small\sffamily\bfseries}%
\def\capt{\small\sffamily}%

\newdimen\@captionmargin\@captionmargin1cm\relax
\newdimen\captionwidth\captionwidth\hsize\relax

\def\eqref#1{(\protect\ref{#1})}

\def\proof{\@ifnextchar[{\@proof}{\@proof[\unskip]}}
\def\@proof[#1]{\noindent{\bf Proof #1.}\enspace}

\def\@mt#1{\ifmmode#1\else$#1$\fi}
\def\qed{\hfill\@mt{\Box}}
\def\qqed{\hfill\@mt{\Box\enspace\Box}}

\let\Bbb\bf

\def\cU{{\cal U}}
\def\cP{{\cal P}}
\def\tg{{\tilde g}}
\def\tZ{{\tilde Z}}
\def\fg{{\frak g}}
\def\tr{\text{tr}}
\def\cZ{{\cal Z}}
\def\cD{{\cal D}}
\def\bR{{\Bbb R}}
\def\bC{{\Bbb C}}
\def\cE{{\cal E}}
\def\bZ{{\Bbb Z}}
\def\bN{{\Bbb N}}
\def\bQ{{\Bbb Q}}
\def\QI{{\Bbb Q}\cup\{\infty\}}
\def\RI{{\Bbb R}\cup\{\infty\}}

\def\bysame{\same[\kern2cm]\,}

\def\br#1{\left\lfloor#1\right\rfloor}
\def\BR#1{\left\lceil#1\right\rceil}

\def\abstractname{}

{\let\@noitemerr\relax
\vskip-2.7em\kern0pt\begin{abstract}
\noindent{\bf Abstract.}\enspace
Introducing a way to modify knots using $n$-trivial rational tangles,
we show that knots with given values of Vassiliev invariants
of bounded degree can have arbitrary unknotting number (extending
a recent result of Ohyama, Taniyama and Yamada). The same result is
shown for 4-genera and finite reductions
of the homology group of the double branched cover. Closer consideration
is given to rational knots, where it is shown that
the number of $n$-trivial rational knots of at most
$k$ crossings is for any $n$ asymptotically at least $C^{(\ln k)^2}$
for any $C<\sqrt[2\ln 2]{e}$.
\\[1mm]
\end{abstract}
}

\section{\label{sect1}Introduction}

In \cite{Stanford}, Stanford introduced a way to modify
knots into alternating prime ones using 3 braids \cite{BirMen},
not affecting (i.~e., changing the values of)
any finite number of Vassiliev invariants
\cite{Birman,BirmanLin,BarNatanVI,BarNatanBibl,%
BarNatanOde,Vassiliev,Vogel}. The 3 braids were chosen to be
iterated pure braid commutators and so they are $n$-trivial in the sense
of Gousarov \cite{Gousarov}, see \cite{bseq}.

In this paper, we give another such construction by means of rational 
tangles, which we describe in section \reference{Sc2}.
It can be applied to any diagram of a knot, not only to closed braid
diagrams. While Stanford's construction is useful not to augment
the braid index (if it is $\ge 3$), our construction
is useful, when applied in an arborescent diagram, not to spoil
arborescency of a knot. Hence, a similar argument to Stanford's
allows us to prove an `arborescent' version of his modification theorem:

\begin{theorem}\label{cr1}
Let $v_1,\dots,v_n$ be Vassiliev invariants. Then for any knot $K$
there is some prime alternating knot $K'$ with $v_i(K)=v_i(K')$ for 
$1\le i\le n$. If $K$ is arborescent, then $K'$ can be chosen to be so
as well.
\end{theorem}

In Gousarov's language two knots $K_1$ and $K_2$ having the same
Vassiliev invariants of degree up to $n$ are called $n$-similar.
We denote this by $K_1\sim_nK_2$.

Applied to the 1-crossing-diagram of the unknot, our method produces
(infinite) series of $n$-trivial 2-bridge knots for given $n\in\bN$.
Hence we have

\begin{corollary}\label{th1}
For any $n$ there exist infinitely many $n$-trivial rational knots
of genus $2^n$. Infinitely many of them have unknotting number one.
\end{corollary}

The number of such knots will be (asymptotically) estimated more
accurately in section \reference{Sc3}. The important feature of this
estimate is that it is asymptotically independent on the degree of
triviality. Such an
estimate does not appear to have been known before (see remark
\reference{rmes}).


Our knots differ in several regards from previous constructions.
Lin's iterated Whitehead doubles \cite{Lin} have genus and unknotting
number one, and are non-alternating, Ng's knots \cite{Ng} are slice
and of unknotting number at most two but their genus is difficult
to control, the same being true for Stanford's alternating braid knots.

Ng's construction offers an analogy to another outcome of our
work. She showed that, beside the Arf invariant, Vassiliev invariants
give no information on knot cobordism. This helps completing a picture,
realized soon after Vassiliev invariants became popular, that all
classical knot invariants (that is, those known before the ``polynomial
fever'' \cite[Preface]{Rolfsen} broke out with \cite{Jones}), are
not, or stronger (almost) unrelated to, Vassiliev invariants, see
\cite{Birman}. Our method exhibits the same picture for the unknotting
number.


\begin{theorem}\label{thM}
Let $K$ be some knot and $n$, $u$ positive integers. Then there
exists a prime knot $K_{n,u}$ of unknotting number $u$ having the same
Vassiliev invariants of degree up to $n$ as $K$. Moreover, for
fixed $K$ and $n$, $K_{n,u}$ can be chosen to be 
alternating (and prime) for almost all $u$.
\end{theorem}

We show this result in section \reference{Sc4}. 
It extends the result of Ohyama-Taniyama-Yamada \cite{OTY}
(see also \cite{Ohyama}), which is the claim of the theorem for $u=1$.
Their result is used in the proof, together with an application of our
method, given $K$, how to construct $K_{n,u}$ for any $u\ge u(K)$.
The use of the tangle calculus of \cite{KL} allows to ensure
primality in most cases, contrarily to Ng's knots, which are composite.
Since we will use the signature for the proof of theorem
\reference{thM}, the same statement holds via the Murasugi--Tristram
inequality also for the 4-genus $g_s>0$ instead of the
unknotting number, thus extending the case $g_s=0$ studied by Ng.


Theorem \reference{thM} is a bit surprising, as the picture changes
when considering other unknotting operations, or, at least
conjecturally, special
classes of knots, see \S\reference{Sc6}. Also, the situation differs
when considering infinitely many Vassiliev invariants, because for
example the Jones polynomial, which by \cite{BirmanLin} is equivalent
to such a collection, does carry some (albeit modest) unknotting number
information, see \cite{LickMil,granny,Traczyk}.


Beside signatures or 4-ball genera, for the unknotting number 
results we use the estimate of Wendt \cite{Wendt}, the number of
torsion coefficients of the homology $H_1(D_K,\bZ)$ of the
double branched cover $D_K$ of $S^3$ over $K$. As a by-product,
we obtain a similar result to the ones above regarding the
homology of $D_K$ over rings of positive characteristic (see theorem
\reference{thpc}). It would be more interesting (but much more
difficult) to examine the situation with the whole $\bZ$-module
$H_1(D_K,\bZ)$.

%

Finally, in section \reference{Sc6}, we conclude by summarizing some
problems suggested by our results.

\section{Rational tangles\label{Sc2}}

In this section we introduce the type of (rational) tangles which will
be applied in the subsequent constructions.

Rational tangles were introduced by Conway \cite{Conway}. The Conway
notation $C(a_1,\dots,a_n)$ of a rational tangle is
a sequence of integers, to which a
canonical diagram of the tangle is associated, see \cite[\S 2.3]{Adams}.
Define the iterated fraction (IF) of a sequence of integers
$a=(a_1,\dots,a_n)$ recursively by
\[
IF(a_1):=a_1\,,\dots,IF(a_1,\dots,a_{n-1},a_n):=\frac{1}{IF(a_1,\dots,
a_{n-1})}+a_n\,.
\]
It will be helpful to extend the operations `$+$' and `$1/.$' to
$\QI$ by $1/0=\infty,\ 1/\infty=0,\ k+\infty=\infty$ for any $k\in
\bQ$. The reader may think of $\infty$ as the fraction $1/0$, to which
one applies the usual rules of fraction arithmetics and reducing.
In particular reducing tells that $-1/0=1/0$, so that for us
$-\infty=\infty$. This may appear at first glance strange, but has
a natural interpretation in the rational tangle context. A rigorous
account on this may be found in Krebes's paper \cite{Krebes}.

In this sense, $IF$ is a map $(\fa n\in \bN)$
\[
IF\,:\,\bZ^n\,\lra\,\QI\,.
\]

It is known \cite{Adams}, that diagrams of sequences of integers
with equal $IF$ belong to the same tangle (up to isotopy; where isotopy
is defined by keeping the endpoints fixed). The correspondence is
\[
C(a_1,\dots,a_n)\,\llra\,IF(a_n,\dots,a_1)\,.
\]

\begin{figure}[htb]
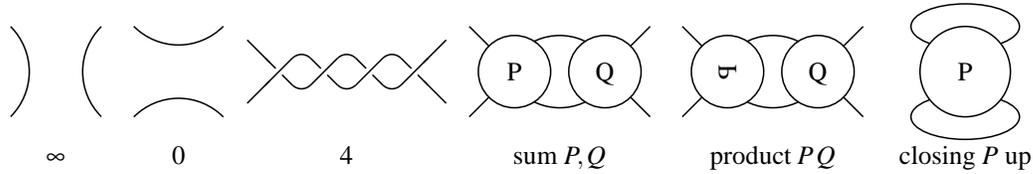

\[
\begin{array}{*6c}
\diag{6mm}{2}{2}{
  \piccirclearc{-1 1}{1.41}{-45 45}
  \piccirclearc{3 1}{1.41}{135 -135}
} & 
\diag{6mm}{2}{2}{
  \piccirclearc{1 -1}{1.41}{45 135}
  \piccirclearc{1 3}{1.41}{-135 -45}
} &
\diag{6mm}{4}{1}{
  \picmultigraphics{4}{1 0}{
    \picline{0.2 0.8}{0.8 0.2}
    \picmultiline{-6 1 -1.0 0}{0.2 0.2}{0.8 0.8}
  }
  \picmultigraphics{3}{1 0}{
    \piccirclearc{1 0.6}{0.28}{45 135}
    \piccirclearc{1 0.4}{0.28}{-135 -45}
  }
  \picline{-0.2 -0.2}{0.3 0.3}
  \picline{-0.2 1.2}{0.3 0.7}
  \picline{4.2 -0.2}{3.7 0.3}
  \picline{4.2 1.2}{3.7 0.7}
} &
\diag{6mm}{4}{2}{
  \picline{0 0}{0.5 0.5}
  \picline{0 2}{0.5 1.5}
  \picline{4 0}{3.5 0.5}
  \picline{4 2}{3.5 1.5}
  \piccirclearc{2 0.5}{1.3}{45 135}
  \piccirclearc{2 1.5}{1.3}{-135 -45}
  \picfilledcircle{1 1}{0.8}{P}
  \picfilledcircle{3 1}{0.8}{Q}
} &
\diag{6mm}{4}{2}{
  \picline{0 0}{0.5 0.5}
  \picline{0 2}{0.5 1.5}
  \picline{4 0}{3.5 0.5}
  \picline{4 2}{3.5 1.5}
  \piccirclearc{2 0.5}{1.3}{45 135}
  \piccirclearc{2 1.5}{1.3}{-135 -45}
  \pictranslate{1 1}{
    \picrotate{90}{
      \picscale{-1 1}{
	\picfilledcircle{0 0}{0.8}{P}
      }
    }
  }
  \picfilledcircle{3 1}{0.8}{Q}
} &
\diag{6mm}{2.4}{3}{
  \picmultigraphics{2}{0 2}{
    \picellipse{1.2 0.5}{1.2 0.5}{}
  }
  \picfilledcircle{1.2 1.5}{1}{P}
} \\
\infty & 0 & 4 & \text{ sum $P,Q$ } &
\text{ product $P\,Q$ } & \text{ closing $P$ up }
\end{array}
\]
\caption{Operations with rational tangles\label{fig1}}
\end{figure}

Using this fact, one can convince himself, that a rational tangle
$T$ has a diagram which closes (in the way described in \cite{Adams},
see also figure \ref{fig1})
to an alternating reduced prime diagram of a link (the only
exception for reducedness being the tangle with notation $(1)$),
which has $1$ or $2$ components (as in our examples below).
This diagram is obtained by taking a representation of $IF(a)=IF(c)$
for a Conway notation $a$ of $T$, such that all numbers in $c$ are of
the same sign (it is easy to see that such a sequence $c$ always
exists). In particular,
$|c|:=\sum_i|c_i|$ is the crossing number of the closure of $T$
(see \cite{Kauffman,Murasugi,Thistle}), and so $T$ is trivial, i.~e.,
the $0$-tangle, iff $IF(a)=0$ (as for $a$ the $0$-tangle and $IF(a)\ne
0$ we had $|c|>0$ and $c\ne (1)$, and thereby a contradiction).
We also see this way, that rational links are prime (see
\cite{Menasco}; this result independently follows from the additivity
of the bridge number proved by Schubert \cite[p.~67]{Adams}).

Define for a finite sequence of integers $a=(a_1,\dots,a_n)$ its
reversion $\ol a:=(a_n,\dots,a_1)$ and its negation by $-a:=(-a_1,\dots,
-a_n)$. For $b=(b_1,\dots,b_m)$ the term $ab$ denotes the concatenation
of both sequences $(a_1,\dots,a_n,b_1,\dots,b_m)$. We also write
$IF(a,a_{n+1})$ for $IF(a_1,\dots,a_n,a_{n+1})$, and analogously
$IF(a,a_{n+1},b)$ etc.

\begin{proposition}\label{th2}
Fix some even $a_1,\dots,a_n\in\bZ$ and build inductively
the integer sequences $w_n$ by
\begin{eqn}\label{was}
w_1:=(a_1),\quad\dots\quad w_n:=w_{n-1}(a_n)\ol{-w_{n-1}}\,.
\end{eqn}
Then the rational tangles with Conway notation $w_n$ are $n$-trivial,
and, if all $a_i\ne 0$, non-trivial, i.~e., not (isotopic to)
the $0$-tangle.
\end{proposition}

\proof For given $n$ consider the braiding polynomial $P$ \cite{bseq}
of some Vassiliev invariants (which may be assumed to be zero on the
$0$-tangle), on the $w_n$-tangle as polynomial in $a_1,\dots,a_n$.
By the discussion of Stanford's examples in \cite{bseq}
and the previous remarks, we need to show that
$IF(\dots,0,\dots)\equiv 0$, and so $P\big|_{a_i=0}\equiv 0$
$\forall i\le n$, and $IF(a_1,\dots,a_n)\ne 0$, if all $a_i\ne 0$.  

Do this by the inductive assumption over $n$. For $n=1$ the claim is evident.
For fixed $n$ by induction assumption $IF\big|_{a_i=0}\equiv 0$
$\forall i< n$, as $IF(a,a_n,\ol{-a})$ is independent of $a_n$ if
$IF(a)=0$, and $IF(a,0,\ol{-a})=0$ for any integer sequence $a$.
But therefore also $IF\big|_{a_n=0}\equiv 0$.

To see that for $a_1\ne 0,\dots,a_n\ne 0$ the tangle is non-trivial,
use that by induction for $a_1,\dots,a_{n-1}\ne 0$ we have
$IF(w_{n-1})\ne 0$ and that therefore the map
\[
a_n\quad\longmapsto\quad IF\left(a_n+\frac{1}{IF(w_{n-1})},
\ol{-w_{n-1}}\right)
\]
is a bijection of $\QI$, so $a_n=0$ can only be a unique zero. \qed

\begin{exam}
For $a_1=2,$ $a_2=-4$ and $a_3=2$ we have $w_1=(2)$, $w_2=(2,-4,-2)$
and $w_3=(2,-4,-2,2,2,4,-2)$.
\end{exam}

\section{Modifying knots\label{Sc3}}

Prepared with the above tangles, we can now describe our modification
technique.

Proposition \reference{th2} already allows to prove the special case of
theorem \reference{cr1} given in the introduction as corollary
\reference{th1}. We first give this proof, before going to prove
theorem \reference{cr1} itself.

\proof[of corollary \reference{th1}]
Corollary \reference{th1} follows directly from the proposition
\reference{th2} by
replacing the $0$-tangle in the unknot diagram $C(0,c)$ for any 
$c\in\bZ$ by some of the tangles in question.
To see that indeed infinitely many
examples arise this way, take $c$ even and use the well-known fact
that the expression of a rational knot with all Conway coefficients
even is unique. The number of even entries is known to be equal
to twice the genus, hence the genus is as asserted. We obtain the
unknotting number property by taking $a_n=\pm 2$. \qed

%
%
We can now prove the arborescent refinement of Stanford's
result from our setting.

\begin{figure}[htb]
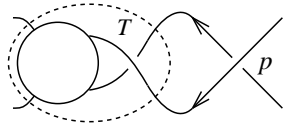

\[
  \diag{6mm}{4}{2}{
    \pictranslate{3 0.0}{
      \picscale{2 2}{
      \picmultivecline{-10 1 -1.0 0}{1 0}{0 1}
      \picmultivecline{-10 1 -1.0 0}{1 1}{0 0}
      }
      \picputtext[l]{1.45 1.0}{$p$}
    }
    \pictranslate{-2 0}{
      \piccurve{5 2}{4.7 2.3}{4.3 2}{4.0 1.5}
      \piccurve{4.0 1.5}{3.7 1}{3.4 0.5}{2.7 0.4}
      \picmulticurve{-10 1 -1.0 0}{4.0 0.5}{3.7 1}{3.4 1.5}{2.7 1.6}
      \picmulticurve{-10 1 -1.0 0}{5 0}{4.7 -0.3}{4.3 0}{4.0 0.5}
      \pictranslate{-1 0}{
        \piccurve{2 0}{2.3 0}{2.3 0}{3 1}
        \piccurve{2 2}{2.3 2}{2.3 2}{3 1}
        \picfilledcircle{3 1}{0.9}{}
	\piclinedash{0.1}{0.05}
	\picellipse{3.7 1}{1.8 1.3}{}
	\picputtext{4.5 1.8}{\small $T$}
      }
    }
}
\]
\caption{Plugging in $T$\label{fig2}}
\end{figure}

\proof[of theorem \reference{cr1}]
Given a knot $K$, take some reduced non-composite diagram  of
$K$ (which exists even if $K$ is composite) and choose a set $S$
of crossings, which need to be switched to obtain from it
an alternating diagram. Then near each such crossing $p$ plug in an
$n$-trivial rational tangle $T$ (in a diagram with alternating
closure), so that the right-most crossing of $T$ cancels with
$p$ by a Reidemeister 2 move (see fig.~\ref{fig2}).

By applying this modification at all crossings in $S$, we are done. \qed

\begin{rem}
More generally, this construction shows that one can preserve
the Conway basic polyhedron.
\end{rem}

We conclude this section with the announced more specific enumeration
result concerning the knots in corollary \reference{th1}.

\begin{corr}
For any $n_o\in\bN$ the number of $n_o$-trivial rational knots
of at most $k$ crossings is asymptotically at least
$C^{(\ln k)^2}$ for any constant $C<\sqrt[2\ln 2]{e}$.
\end{corr}

\proof First note, that the freedom to vary $C$ allows us to replace
for convenience $k$ by $k/2$, or equivalently to consider at most
$2k$ crossing diagrams (instead of at most $k$ crossings).

A diagram of the kind constructed in the proof of proposition
\reference{th2} with $2k$ crossings in the groups of twists
except the first one corresponds to writing
\[
k\,=\,\sum_{i=0}^n\,2^i|w_i|
\]
for some $w_i\in\bZ$ and $n\in\bN$. For $n\le n_o$ the number of
such representations in polynomially bounded in $k$, hence, assuming
we can show the lower bound for the diagrams including these with $n\le
n_o$, it
is possible to neglect them and assume $n\ge n_o$, so that all diagrams
are $n_o$-trivial.

Let
\[
D_k\,:=\,\left\{\,(w_0,\dots,w_n)\,:\,k=\sum_{i=0}^n\,2^i|w_i|,\,
w_i\ne 0,\, n>0\,\right\}
\]
and $d_k:=\#D_k$. Then $d_1=2$ and
\[
d_k\,=\,2\sum_{i=1}^{\br{k/2}} d_i\mbox{\quad for $k\ge 2$.}
\]
To prove the corollary it suffices to show that 
\begin{eqn}\label{1*}
c_k\,:=\,C^{(\ln k)^2}\,<\,(2-\eps)\,\sum_{i=1}^{\br{k/2}} C^{(\ln i)^2}
+D
\end{eqn}
for some $D\in\bR$, $\eps>0$ and sufficiently large $k$, as then
(for possibly larger $k$) $D<C^{(\ln k)^2}\cdot\frac{\eps}{2}$, so
\[
C^{(\ln k)^2}\,<\,2\,\sum_{i=1}^{\br{k/2}} C^{(\ln i)^2}\,,
\]
and hence $d_k\ge C'\cdot c_k+C''$ (for some $C',C''\in\bR$, $C'>0$),
but $C'$ and $C''$ can be eliminated by varying $C$.

To show \eqref{1*}, first use that $i\mapsto C^{(\ln i)^2}$ is 
monotonously growing for $i\ge 1$, so 
\begin{eqn}\label{3*}
\int_{1}^{(k-1)/2}\,C^{(\ln t)^2} dt\,<\,
\sum_{i=1}^{\br{k/2}} C^{(\ln i)^2}\,.
\end{eqn}
Now for \eqref{1*} it suffices to show the inequality for the
derivations of the left hand-sides of \eqref{1*} and \eqref{3*}
for sufficiently large $k$.

But putting $C=e^p$ with $p<\frac{1}{2\ln 2}$, we have that
\[
\frac{d}{dk}\left(\,e^{p(\ln k)^2}\right)\,<\, (2-\eps)\,
\frac{d}{dk}\left(\int_{1}^{(k-1)/2}\,e^{p(\ln t)^2} dt\,\right)
\]
is equivalent to
\[
e^{p(\ln k)^2}\frac{2\ln k}{k}p\,<\,
\frac{2-\eps}{2}e^{p(\ln(k-1)-\ln 2)^2}\,,
\]
and logarithming we get
\[
p(\ln k)^2+\ln 2+\ln\ln k-\ln k+\ln p\,<\,
\ln(2-\eps)-\ln 2+p(\ln(k-1))^2-2p\ln 2\ln(k-1)+p(\ln 2)^2\,.
\]
This is for some $D'\in\bR$ the inequality
\[
p\left((\ln k)^2-(\ln(k-1))^2\right)+\ln\ln k\,<\,
(1-2p\ln 2)\ln(k-1)+D'\,.
\]
Now as $k\to\infty$, the first term on the left goes to 0,
and then the claim is obvious from the condition on $p$. \qed

\begin{rem}\label{rmes}
It should be remarked that the asymptotical estimate itself does
\em{not} depend on $n$. Such an unconditional statement does not
seem to have been known before. For example, the number of Lin's
iterated $n$-fold Whitehead doubles for fixed $n$ grows exponentially
in $k$, because of the result of \cite{Welsh} and the uniqueness of
the companion, but the base of this exponential heavily depends on
$n$-- roughly augmenting $n$ by $1$ requires to take the fourth root
of the base. On the other hand, the dependence on $n$ in our
estimate is present, namely in how quickly the numbers attain their
asymptotical behaviour. Thus our result does not imply the existence
of knots which are $n$-trivial for all $n$. In fact, as our knots
are alternating, no one of them can have this property.
\end{rem}

\section{Unknotting numbers and $n$-triviality\label{Sc4}}

Here we record some consequences of the preceding results concerning
unknotting numbers. The first one is rather easy, and will be
later refined to give a proof of theorem \reference{thM}.

\begin{prop}\label{thu0}
Let $K$ be some knot. Then for any $n\in\bN$ and $u_0\ge u(K)$ there
exists a knot $K_{n,u_0}$ with $u(K_{n,u_0})=u_0$ and $v(K_{n,u_0})=
v(K)$ for any Vassiliev invariant $v$ of degree up to $n$.
\end{prop}

\proof Consider the knots $K_{(i)}:=K\#\bigl(\#^iK'\bigr)$ with
$K'$ being an $n+1$-trivial rational knot of unknotting number one
(provided by corollary \reference{th1}).
Then the Vassiliev invariants of degree up to $n$ of all $K_{(i)}$
are the same as those of $K$, and that any $u_0\ge u(K)$ is the
unknotting number of some $K_{(i)}$ follows from the obvious inequality
$u(K_{(i+1)})\le u(K_{(i)})+1$ and the reverse estimate $u(K_{(i)})\ge
d_{K_{(i)}}=d_K+i$, where $d_K=\rk H_1(D_{K},\bZ)$ is the number of
torsion coefficients of
$H_1(D_K,\bZ)$ and $D_K$ is the double cover of $S^3$ branched along
$K$, see \cite{Wendt}. \qed

%
%
%

Now we indicate how to modify the proof of proposition \reference{thu0}
to signatures and 4-genera. (This can also
be deduced from Ng's work, but the proof is now brief, so we
can give it in passing by.)

\begin{theorem}\label{thgg}
Let $n\in\bN$ and $K$ be some knot. Then
{\nopagebreak
\def\labelenumi{\roman{enumi})}\mbox{}\\[-18pt]
\def\theenumi{\roman{enumi}}

\begin{enumerate}
\item\label{item1} for any $s\in2\bZ$ there is a knot $K_{n,s}\sim_nK$
with $\sg(K_{n,s})=s$.
\item\label{item2} for any integer $g\ge 0$ there is a knot
$K_{n,g}\sim_nK$ with $g_s(K_{n,g})=g$, except if $Arf(K)=1$,
$g=0$ (and $n>1$).
\end{enumerate}
}
\end{theorem}

\proof[sketch] By the result of Ohyama, Taniyama and Yamada, Ng's
work for $g_s=0$ (which we cite, but do not cover with our arguments)
and the previous arguments, together with the standard inequalities
$|\sg/2|\le g_s\le u$ (see \cite{Murasugi2,Tristram}), we
are basically left with showing that
the $n$-trivial rational knot of corollary \reference{th1} can be chosen
to be of signature $\pm 2$. For this we remark that the determinant
shows that the signature of a rational knot $S(p,q)$ with $p,q>0$
(in Schubert's notation \cite{Schubert}) is divisible by $4$ exactly
if $p\equiv 1\bmod 4$. Violating this property reduces to
making the number $c$ in the proof of corollary \reference{th1} small
or large enough in order to adjust the desired sign of $IF(w_n,c)$. \qed

Next, we state and prove the result for
the homology of the double branched cover.

\begin{theorem}\label{thpc}
Let $p>2$ be an odd integer, $H$ be a finite $\bZ_p$-module, $n>0$
be a natural number and $K$ be any knot. Then there is a knot $K'$,
which can be chosen to be prime and alternating, with $K'\sim_n K$,
such that $H_1(D_{K'},\bZ_p)=H$.
\end{theorem}

We start the proof by two lemmas, lemma \reference{lM1}
and \reference{lM2}. Then we prove theorem
\reference{thpc} by taking the connected sum of the knots
constructed in the lemmas. 

\begin{lem}\label{lM1}
Let $p>1$ be an odd integer. Then for any $k\in\bZ_p$ and any $n$
there is an $n$-trivial rational knot $K$ with $\dt(K)\equiv k\bmod p$
and $u(K)=1$.
\end{lem}

\proof Let $w_n$ be the sequences of integers as in \eqref{was}
with all $a_i\in\{\pm 2\}$. Then, by the calculation used in
\cite[proof of theorem 1, (ii)$\Ra$ (iii)]{KanMur} we find
\[
\big|\,IF(w_n)\,\big|\,=\,\frac{2^{2^n-1}}{\sum_{i=1}^{n}\pm
2^{2^n-2^i}}\,
\]
for certain signs in the sum depending on the signs of the $a_i$.

Then for the rational tangle with Conway notation $(w_n,s)$ for
a natural (not necessarily even) number $s$ we have
\[
\dt\bigl(\ol{(w_n,s)}\bigr)=\sum_{i=1}^{n}\pm 2^{2^n-2^i}+
s\cdot 2^{2^n-1}\,,
\]
and the existence of proper choice of $s$ follows from the
fact that $2^{2^n-1}$ and $p$ are relatively prime. \qed

\begin{lem}\label{lmphi}
The (tautological) homomorphism $h_{p,q}:\,\bZ_p^*\to\bZ_q^*$
for any $q\mid p$ is onto (where $\bZ_p^*$ is the group of units
of $\bZ_p$ or the relatively prime to $p$ rest classes modulo $p$).
\end{lem}

\proof This surjectivity follows because $|\bZ_p^*|=\phi(p)$,
$\bZ_{ab}^*=\bZ_{a}^*\times \bZ_{b}^*$ for $(a,b)=1$ and because
obviously $\big|\ker(h_{p^u,p^{u-1}})\big|=p=\phi(p^u)/\phi(p^{u-1})$
for any prime $p$ and $u>1$. \qed

\begin{lem}\label{lM2}
Let $p>1$ be an odd integer and $K$ be an unknotting number one
knot. Then for any $n$ there is a knot $K'\sim_n K$ with
$\dt(K')$ relatively prime to $p$.
\end{lem}

\proof We use the tangle calculus of Krebes (see \cite{Krebes} for
details). He showed that the pair of determinants of the closures of a
tangle $T$ and its flipped version $T0$ (the product of $T$ with the
$0$ tangle in the notation of figure \reference{fig1})
can be viewed as the numerator and denominator of a certain
generalized rational number, denoted here by $R(T)$, lying in
$\tl\bQ:=\bZ\times\bZ \big/(a,b)\sim(-a,-b)$, which (up to signs)
is additive under tangle sum (as in figure \reference{fig1}),
and generalizes $IF$ for rational tangles. 

The fact that $u(K)=1$ shows that $K$ can be presented as the closure
$\ol{T}$ of a tangle $T$ such that the closure $\ol{T,2}$ of the
tangle sum of $T$ with the $2$-tangle (clasp) is the unknot.
Krebes's calculus then shows that $R(T)=(\pm 2k\pm 1)/k\in\tl\bQ$
for certain signs and a natural number $k$.

Then consider the tangle sum of $T$ with the rational tangle $S=(w_n,s,
0)$ for a (not necessarily even) integer $s$.
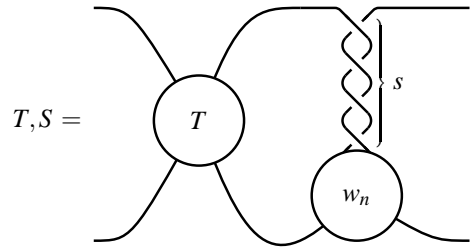
\[
T,S\,=\,
\diag{1cm}{5}{5}{
 \pictranslate{3.75 2}{
  \picrotate{90}{
    \picscale{0.5 d}{
  \picmultigraphics{4}{1 0}{
    \picline{0.2 0.8}{0.8 0.2}
    \picmultiline{-6 1 -1.0 0}{0.2 0.2}{0.8 0.8}
  }
  \picmultigraphics{3}{1 0}{
    \piccirclearc{1 0.6}{0.28}{45 135}
    \piccirclearc{1 0.4}{0.28}{-135 -45}
  }
  \picarcangle{0.3 0.3}{0 0}{-0.4 -0.2}{0.2}
  \picarcangle{0.3 0.7}{0 1}{-0.4 1.2}{0.2}
  \picarcangle{3.7 0.3}{4 0}{4 -2.5}{0.2}
  \picarcangle{3.7 0.7}{4 1}{4 2}{0.2}
  \picputtext{2 0}{$\underbrace{\rx{3.4\unitlength}}_{}$}
}}}  
  \picputtext{4.05 3}{$s$}
  \piccurve{3.5 1.5}{2.5 0.5}{2.1 0.5}{1.4 2.5}
  \piccurve{3.5 1.5}{4.5 0.9}{4.5 0.9}{5 0.9}
  \piccurve{1.4 2.5}{0.5 0.9}{0.5 0.9}{0 0.9}
  \piccurve{1.4 2.5}{0.5 4}{0.5 4}{0 4}
  \piccurve{1.4 2.5}{1.9 4}{2.2 4}{2.75 4}
  \picfilledcircle{1.4 2.5}{0.6}{$T$}
  \picfilledcircle{3.5 1.5}{0.6}{$w_n$}
}
\]
Then $\ol{T,S}\sim_n K$ and by Krebes's calculus 
\[
\dt(\ol{T,S})=k2^{2^n-1}+(\pm 2k\pm 1)\Bigl[\,\sum_{i=1}^{n}
\pm 2^{2^n-2^i}+s\cdot 2^{2^n-1}\,\Bigr]\,.
\]

Changing $s$ by $\pm 1$ causes the expression to change by
$(\pm 2k\pm 1)2^{2^n-1}$. Thus we could finish the proof as in
the case of lemma \reference{lM1} unless $2k\pm 1$ and $p$ are not
relatively prime. In this case let $l=(2k\pm 1,p)$ be their greatest
common divisor. Clearly $(l,k2^{2^n-1})=1$ and so $(l,k')=1$ for
\[
k'=k2^{2^n-1}+(\pm 2k\pm 1)\,\sum_{i=1}^{n}\pm 2^{2^n-2^i}\,.
\]
We would be done if we can find an $s'\in\bZ$ with
$(k'+l\cdot s',p)=1$. Then set 
\[
s:=s'\cdot \frac{1}{2^{2^n-1}} \cdot \frac{l}{2k\pm 1}
\]
in $\bZ_p$. Here the meaning of the second factor is clear, as
$2^{2^n-1}$ is invertible in $\bZ_p$. The third factor means some
(fixed) preimage under $h_{p,p/l}$ of the (multiplicative) inverse
of $(2k\pm 1)/l\in\bZ_{p/l}^*$. The existence of this preimage
follows from lemma \reference{lmphi}. In turn,
the existence of $s'$ is equivalent to the surjectivity of the
homomorphism $h_{p,l}$, which again follows from lemma
\reference{lmphi}. \qed

\proof[of theorem \reference{thpc}]
We can write
\[
H\,=\,\bigoplus_{i=1}^{l}\,\bZ_{p_i}
\]
with $p_i\mid p$. Let $\hat K$ be the knot found to $K$ in lemma
\reference{lM2}, and $K_i$ be the knots from lemma \reference{lM1}
for $k=p/p_i$. Then $H_1(D_{\hat K},\bZ_p)=1$, and $H_1(D_{K_i},\bZ_p)=
\bZ_{p_i}$, since $K_i$ are rational by construction.
Thus
\[
K''\,=\,\hat K\es\#\es\mathop{\#}_{i=1}^l\,K_i
\]
is a knot with the desired values of Vassiliev invariants and
homology group. It remains to make $K''$ into a prime alternating knot
$K'$, which will be the knot we sought.  

To obtain $K'$ from $K''$, take a prime diagram of $K''$, and apply
the plugging technique in the proof of theorem \reference{cr1}
with a tangle $w_n$ of the form \eqref{was} with $a_n=2p$.
Then by the work of Gordon and Litherland \cite{GorLit} on the Goeritz
matrix, this plugging preserves the structure of $H_1(D_{K''},\bZ_p)$,
since $w_n$ turns into the $0$-tangle by changing some of its
Conway coefficients by a multiple of $p$. \qed


\begin{rem}
It is uninteresting to consider $p$ to be even, because for any
\em{knot} (although not link) $K$, $H_1(D_K,\bZ)$ has no 2-torsion,
so its reduction modulo $2p$ is equivalent to its reduction modulo $p$.
\end{rem}

\begin{rem}
Instead of making $K'$ in theorem \reference{thpc} alternating
and prime, we can also achieve, setting $K'=K''$, that it has
$u(K')\le \rk
H_1(D_{K'},\bZ_p)+2$, as the knot in lemma \reference{lM1} had
unknotting number 1, and this constructed in lemma \reference{lM2}
has unknotting number 1 or 2.
\end{rem}

We conclude this section with the proof of theorem \reference{thM}.
For this we use the prime tangle calculus of \cite{KL}. Recall that
a tangle is called prime if it contains no
properly embedded separating disk, and no one of the strands
has a connected summand (i.e. a sphere intersecting it in a
knotted arc). First we need a simple lemma.

\begin{lem}
There are prime tangles with unknotted closure.
\end{lem}

\proof Consider the knot $9_{34}$, which has unknotting number 1,
and the encircled crossing, whose switch unknots it.
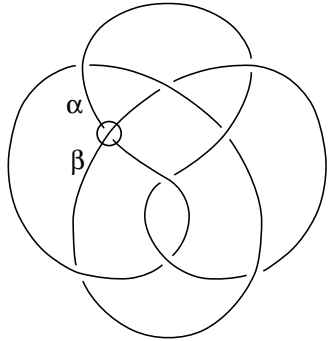
\[
\diag{2cm}{2.2}{2.5}{
  \picputtext[dl]{0 0}{
     \epsfs{4.3cm}{t1-9-34}
  }
  \piclinewidth{90}
  \piccircle{0.73 1.42}{0.08}{}
  \picputtext{0.5 1.6}{$\ap$}
  \picputtext{0.52 1.25}{$\bt$}
}
\]
Switching the crossing, and cutting the edges $\ap$ and $\bt$ we
obtain (up to change of the unbounded region) a tangle $U$
with unknotted closure. To show primeness, we need to show
first that it has no connected summand. However, this is
clear since the closure is unknotted. Then, we need to ensure
that it is not rational. For this consider the other
closure $\ol{U\cdot 0}$ of $U$. It has an alternating diagram
with Conway polyhedron \cite{Conway} $6^*$, and hence it is not
rational. Thus $U$ is not a rational tangle, and is therefore prime.
\qed

\proof[of theorem \reference{thM}]
Fix $K$ and $n$. Let $K_{n,1}$ be the knot constructed in
\cite{Ohyama}. Since $u(K_{n,1})=1$, by \cite{Scharl,Zhang}, $K_{n,1}$
is prime, and thus by \cite{KL}, $K_{n,1}=\ol{T_{n,1}}$, with $T_{n,1}$
being a prime tangle. 
We can without loss of generality assume that the orientation of $T_{n,1}$ is like
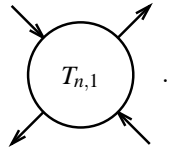
\[
\diag{1cm}{2}{2}{
 \pictranslate{1 1}{
  \picvecline{45 0.7 x polar}{45 1.3 x polar}
  \picvecline{225 0.7 x polar}{225 1.3 x polar}
  \picvecline{135 1.3 x polar}{135 0.7 x polar}
  \picvecline{315 1.3 x polar}{315 0.7 x polar}
  \picfilledcircle{0 0}{0.7}{$T_{n,1}$}
 }
}\,.
\]
Otherwise, we can replace $T_{n,1}$ by its sum with a one-crossing
tangle. This sum is again a prime tangle (see \cite{Van}). Let $w_{n+1
}$ be a $(n+1)$-trivial rational tangle, and $T_n'=w_{n+1}\cdot c_n$.
Let $U$ be a prime tangle with unknotted closure and set $T_n''=
U\cdot T_n'$. Then $T_n''$ is also prime.

Since smoothing out a crossing in the group of $c_n$
gives the link $\ol{w_{n+1}}$, which has non-zero determinant,
as in the proof of theorem \reference{thgg}, 
by choosing $c_n$ large or small enough, we can achieve
that $\sg(\ol{T_n''})\ne 0$. Also, by choosing $a_{n+1}=\pm 2$, we
can achieve that $u(\ol{T_n''})=1$.

Now consider
\[
T_{n,k}\,=\,T_{n,1},\kern0.2ex\underbrace{\kern-0.2ex T_n''\cdot
0, T_n''\cdot 0,\dots,T_n''\cdot 0\kern-0.2ex}_{\scbox{$k-1$ times}}
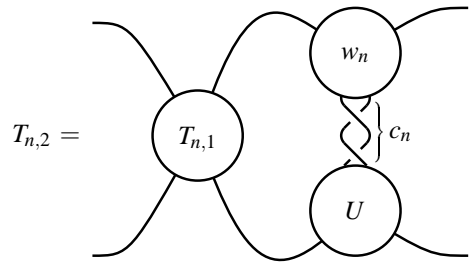
\kern0.2ex\,.
\]
\[
T_{n,2}\,=\,
\diag{1cm}{5}{5}{
 \pictranslate{3.75 2}{
  \picrotate{90}{
    \picscale{0.5 d}{
  \picmultigraphics{4}{1 0}{
    \picline{0.2 0.8}{0.8 0.2}
    \picmultiline{-6 1 -1.0 0}{0.2 0.2}{0.8 0.8}
  }
  \picmultigraphics{3}{1 0}{
    \piccirclearc{1 0.6}{0.28}{45 135}
    \piccirclearc{1 0.4}{0.28}{-135 -45}
  }
  \picputtext{1.1 0}{$\underbrace{\rx{1.6\unitlength}}_{}$}
}}}  
  \picputtext{4.1 2.55}{$c_n$}
  \piccurve{3.5 1.5}{2.5 0.5}{2.1 0.5}{1.4 2.5}
  \piccurve{3.5 1.5}{4.5 0.9}{4.5 0.9}{5 0.9}
  \piccurve{1.4 2.5}{0.5 0.9}{0.5 0.9}{0 0.9}
  \piccurve{1.4 2.5}{0.5 4}{0.5 4}{0 4}
  \piccurve{1.4 2.5}{2.0 4.4}{2.4 4.3}{3.0 3.9}
  \piccurve{3.5 3.6}{4.5 4.2}{4.5 4.2}{5 4.2}
  \picfilledcircle{1.4 2.5}{0.6}{$T_{n,1}$}
  \picfilledcircle{3.5 1.5}{0.6}{$U$}
  \picfilledcircle{3.5 3.6}{0.6}{$w_n$}
}
\] 
We have that 
\begin{eqn}\label{U1}
u(\ol{T_{n,k+1}})\le u(\ol{T_{n,k}})+1\,.
\end{eqn}
Then, because of the above choice of orientation of $T_{n,1}$,
the tangle $T_{n,k}$ differs from $T_{n,k},\infty$ by a band
connecting (plumbing of a Hopf band). But the closure of
$T_{n,k},\infty$ is $\ol{T_{n,1}\cdot 0}\,\#\,(\#^{k-1}\ol{T_n''})$,
and since $\sg(\ol{T_n''})\ne 0$, we have
\[
2u(\ol{T_{n,k}})\,\ge\,|\,\sg(\ol{T_{n,k}})\,|\,\ge\,
\big|\,\sg\bigl(\ol{T_{n,1}\cdot 0}\#\,(\#^{k-1}\ol{T_n''})\bigr)\,\big|
-1\,\ge\,(k-1)|\,\sg(\ol{T_n''})\,|-|\,\sg(\ol{T_{n,1}\cdot 0})\,|-1
\,\lra\,\infty\,,
\]
as $k\to\infty$. This, together with \eqref{U1} and $u(K_{n,1})=1$,
shows that each natural
number $u$ is realized as the unknotting number of some
$\ol{T_{n,k}}$, with $k\ge u$. Since
\[
\kern0.2ex\underbrace{\kern-0.2ex T_n''\cdot 0,
T_n''\cdot 0,\dots,T_n''\cdot 0\kern-0.2ex}_{\scbox{$k$ times}}
\kern0.2ex
\]
is prime for $k\ge 1$ by \cite{Van}, $\ol{T_{n,k}}$ is a prime knot
for $k\ge 2$, and also for $k=1$ by \cite{Scharl}.

To show the claim for prime alternating knots, it suffices to
replace in the above argument $T_{n,1}$ by an alternating
tangle $\hat T_{n,1}$, obtained from $T_{n,1}$ by the
operation described in the proof of theorem \reference{cr1}
(and on figure \reference{fig2}), and to take instead of
$T_n''$ the alternating tangles $T_n'$, mirrored
in such a way so as $T_{n,k}$ to remain alternating. \qed

\section{Odds \& Ends\label{Sc6}}

There are a lot of questions and problems suggested by the
above results. Here we give an extensive summary of what one
could think about to improve and push further.

We start by a problem concerning the construction itself.

\begin{question}
Although they easily achieve alternation, both our and Stanford's 
constructions live at the cost of exponential (in $n$) crossing number
augmentation (at least in the diagrams where $n$-triviality
is achieved). Contrarily, the series of examples of $n$-trivial
knots by Ng \cite{Ng} have crossing number which is
linearly bounded in $n$. There knots are, however, not (\em{a priori})
alternating or positive, slice (so all have zero signature),
and so not to distinguish among each other by such \em{ad hoc}
arguments as \cite{Kauffman,Murasugi,Thistle}. Is it possible to
combine the advantages of both series of examples in a new one?
\end{question}

As for the applications of our construction,
the results of the section \reference{Sc4} suggest
two more problems.

\begin{question}
Does an alternating prime knot $K_{n,u_0}$ exist for any choice of
$n$ and $u_0$ in theorem \reference{thu0}?
\end{question}

Theorem \reference{thM} can be interpreted as saying that any finite
number of Vassiliev invariants does not obstruct to any (non-zero)
value of the unknotting number. On the other hand,
%
it is remarkable that such obstructions do exist for
other unknotting operations, as the $\Dl$ move of Murakami and
Nakanishi \cite{MurNak}. Moreover, certain properties of
Vassiliev invariants with respect to the ordinary unknotting operation
can be suspected in special cases.

\begin{conjecture}
Let $(K_i)$ be a sequence of (pairwise distinct)
positive knots of given unknotting number,
and $v_2=-\myfrac{1}{6}V''(1)=\myfrac{1}{2}\Dl''(1)$
and $v_3=-\myfrac{1}{12}V''(1)-\myfrac{1}{36}V'''(1)$
be the (standardly normalized) Vassiliev invariants
of degree $2$ and $3$, where $V$ is the Jones \cite{Jones} and
$\Dl$ the Alexander polynomial \cite{Alexander}. Then the numbers
$\log_{v_2(K_i)}v_3(K_i)$ (which are well-defined for $K_i\ne !3_1$)
converge to $2$ as $i\to\infty$.
\end{conjecture}

This conjecture is related to some results of \cite{gen1}, but
it would take us too far from the spirit of this paper to describe 
the relation closer here.

\begin{question}
Can any finite number of Vassiliev invariants be realized by a
quasipositive (or strongly quasipositive) knot?
\end{question}

\begin{rem}
Rudolph showed that any Seifert pairing can be realized by
a quasipositive knot, so there are no constraints to
quasipositivity from Vassiliev invariants via the Alexander polynomial.
\end{rem}

The consideration of the  homology of the double branched cover
suggests several questions about further generalizations and
modifications, basically coming from the desire to remove the
reduction modulo some number. We should remark that $n$-similarity poses
via the Alexander polynomial a congruence condition on the determinant,
and ask whether this is the only one.

\begin{question}
Is there for any $n\in\bN$ and any knot $K$ a knot $K'$ with
$K'\sim_nK$ (or weaker, an $n$-trivial knot $K'$) with
{\nopagebreak
\def\labelenumi{\roman{enumi})}\mbox{}\\[-18pt]
\def\theenumi{\roman{enumi}}

\begin{enumerate}
\item\label{Item1} any (finite) abelian group of order $\dt(K')
\equiv \pm\dt(K)\bmod 4^{\br{(n+1)/2}}$ as
homology of the double branched cover, or weaker
\item\label{Item2} any odd positive integer 
$\dt(K')\equiv \pm\dt(K)\bmod 4^{\br{(n+1)/2}}$ as determinant?
\end{enumerate}
}
As a weaker version of part \reference{Item1}), is any (non-constant)
knot invariant depending (only) on the homology of the double branched
cover not a Vassiliev invariant?
\end{question}

\begin{rem}
Note, that there is no chance to get $K'$ with some of the
above properties in general to be alternating, as for $K$ alternating
$\dt(K)\ge c(K)$.
\end{rem}

\begin{rem}
The weaker statement that any non-constant knot invariant depending
on finite reductions of the homology of the double branched
cover is not a Vassiliev invariant is, as seen, true, and beside from 
its stronger versions proved above, originally follows
from the results on $k$-moves in \cite{bseq}. 
\end{rem}

At least for part \reference{Item2}) the strategy
followed in \S\reference{Sc4} appears promising~-- use Krebes calculus
and construct arborescent tangles by properly inserting $n$-trivial
rational tangles. This leads to a question on the image of Krebes's
invariant $R$ on the set of $n$-trivial arborescent tangles, whose first
part is a specialization of \reference{Item2}) of the question above,
and whose second part addresses another unrelated by appealing property.

\begin{question}
Let $T_n$ be the set of $n$-trivial arborescent tangles
(of the homotopy type of the $0$-tangle).
{\nopagebreak
\def\labelenumi{\roman{enumi})}\mbox{}\\[-18pt]
\def\theenumi{\roman{enumi}}

\begin{enumerate}
\item\label{ITem1} Can for any $n$ and any odd $c$ (or weaker $c=1$)
some $(d,c)\in\tl\bQ$ be realized as $R(T)$ for some $T\in T_n$?
\item\label{ITem2} Is the image of $R(T_n)$ under the 
(tautological) homomorphism $\tl\bQ\to\QI$ dense in $\RI$? Is it even
the whole $\QI$?
\end{enumerate}
}
\end{question}

\begin{rem}
Note, that Krebes in his paper (basically) answers positively both
questions in part \reference{ITem2}) for $n=0$.
\end{rem}

%
%

In view of the desire to consider the detection of orientation
(which is a much more relevant problem than just the detection of
knottedness), the constructions of $n$-similar knots
suggest one more general and final problem.

\begin{question}
In \cite{bseq}, I gave a generalization of Gousarov's concept of
$n$-triviality, called $n$-in\-ver\-tibi\-lity, which, \em{inter alia},
led by use of \cite{BirMen} to an elementary construction of a 14
crossing (closed) 3-braid knot, whose orientation cannot be detected
in degree $\le 11$, the argument being provided without any computer
calculation. The argument applied there does not seem (at least
straightforwardly) to be recoverable from $n$-triviality alone (in
particular, because the knot is not $11$-trivial). However, yet,
I have no series of examples of arbitrary degree (as those here),
where this generalized argumentation shows indeed more powerful,
that is, where the failure of Vassiliev invariants of degree $\le n$
to detect orientation can be explained via $n$-invertibility, but not
via $n$-similarity to some invertible knot. Do such examples exist?
\end{question}

\noindent{\bf Acknowledgement.} I would wish to thank to N.~Askitas
for interesting discussions and helpful remarks on unknotting numbers,
and to P.~Traczyk for sending me a copy of \cite{Bernhard}.

\end{document}